\documentclass[a4paper,11pt]{amsart}

\usepackage{amsthm}
\usepackage{amsmath}
\usepackage{amssymb}
\usepackage{latexsym}
\usepackage{enumerate}

\usepackage{graphicx}

\renewcommand{\thetheoremName}

\newtheorem{theorem}{Theorem}[section]
\newtheorem{lemma}[theorem]{Lemma}

\newtheorem{proposition}[theorem]{Proposition}
\newtheorem{corollary}[theorem]{Corollary}

\theoremstyle{definition}
\newtheorem{definition}[theorem]{Definition}

\newtheorem{remark}[theorem]{Remark}

\numberwithin{equation}{section}

\newcommand{\D}{\operatorname{D}}
\newcommand{\Hess}{\operatorname{Hess}}
\newcommand{\dist}{\operatorname{dist}}
\newcommand{\Vol}{\operatorname{Vol}}

\newcommand{\Div}{\operatorname{div}}
\newcommand{\C}{\operatorname{Cap}}

\newcommand{\erre}{\mathbb{R}}
\newcommand{\Lmod}{\operatorname{L}}

\newcommand{\pL}{{\operatorname{\Delta}}_{p}}

\newcommand{\pLP}{{\operatorname{\Delta}}_{p}^{S}}

\newcommand{\pCap}{{\operatorname{Cap}}_{p}}
\newcommand{\loc}{{\rm loc}}
\begin{document}

\title[$p$-Capacity and $p$-Hyperbolicity]{$p$-Capacity and $p$-Hyperbolicity\\
of Submanifolds}

\author[I. Holopainen]{Ilkka Holopainen$^{\dag}$}
\address{P.O. Box 68 (Gustaf H\"allstr{\"o}min katu 2b) \\
FIN-00014 University of Helsinki, Finland.}
\email{ilkka.holopainen@helsinki.fi}
\thanks{$^{\dag}$Work partially supported by the Academy of Finland, project 53292.}

\author[S. Markvorsen]{Steen Markvorsen$^{\#}$}
\address{Department of Mathematics, Technical University of Denmark.}
\email{S.Markvorsen@mat.dtu.dk}
\thanks{$^{\#}$Work partially supported by
the Danish Natural Science Research Council and DGI grant
MTM2004-06015-C02-02.}

\author[V. Palmer]{Vicente Palmer*}
\address{Departament de Matem\`{a}tiques, Universitat Jaume I, Castellon,
Spain.} \email{palmer@mat.uji.es}
\thanks{* Work partially supported by
the Caixa Castell\'{o} Foundation, DGI grant MTM2004-06015-C02-02,
and by the Danish Natural Science Research Council.}



\maketitle
\bibliographystyle{acm}
\begin{abstract}
We use drifted Brownian motion in warped product
model spaces as comparison constructions to show
$p$-hyperbolicity of a large class of submanifolds
for $p\ge 2$.
The condition for $p$-hyperbolicity is expressed in
terms of upper support functions for the radial
sectional curvatures of the ambient space and for
the radial convexity of the submanifold. In the
process of showing $p$-hyperbolicity we also obtain
explicit lower bounds on the $p$-capacity of
finite annular domains of the submanifolds in
terms of the drifted $2$-capacity of the
corresponding annuli in the respective comparison
spaces.
\end{abstract}

\section{Introduction} \label{secIntro}
In \cite{Hdir} the first named author solved the asymptotic Dirichlet
problem  at infinity for the $p$-Laplacian in Cartan--Hadamard
manifolds of pinched negative sectional curvature. As a consequence,
such a manifold admits a wealth of non-constant bounded $p$-harmonic
functions. On the other hand, there are no non-constant positive
$p$-harmonic functions on a complete \linebreak Riemannian manifold with
non-negative Ricci curvature; see e.g. \cite{CHS}.
The purpose of the present paper is to initiate the study
of the $p$-Laplacian and the existence of $p$-harmonic functions of various types
on submanifolds. In this paper we concentrate on $p$-hyperbolicity of
submanifolds.

To describe the problem
we are dealing with, suppose that $S$ is a Riemannian
submanifold of an
ambient Riemannian manifold $N$. We look for the most general intrinsic
geometric condition on $N$ and the most general extrinsic geometric
condition on $S$ which together will assure that $S$ is $p$-hyperbolic.
Recall that a Riemannian manifold $M$ is called {\em{$p$-hyperbolic}}, with
$1<p<\infty$, if there exists a compact set $K\subset M$ of positive
$p$-capacity $\pCap (K,M)$ relative to $M$. Here the $p$-capacity of $K$
is defined by
\[
\pCap (K,M)=\inf_{u}\int_{M}\Vert \nabla u\Vert^{p}\,d\mu,
\]
where the infimum is taken over all real-valued
functions $u\in C^{\infty}_{0}(M)$, with $u\ge 1$
in $K$. In case $p=2$, the $p$-hyperbolicity of
$M$ is equivalent both to the existence of a
positive Green's kernel for the Laplace-Beltrami
operator and to the transience of $M$, (see the
works \cite{MR755228} and \cite{Gri}). Using the
particular 2-capacity condition alluded to above,
the two last named authors have obtained
geometric criteria for $2$-hyperbolicity of
minimal - or close to minimal - submanifolds in
manifolds with sectional curvatures bounded from
above, (see \cite{MP2} and \cite{MP3}).

In the general case of $1<p<\infty$, the
$p$-hyperbolicity of $M$ is known to be
equivalent to the existence of a (positive)
Green's function $g=g(\cdot,y)$ for the
$p$-Laplace equation, i.e. a certain positive
solution (in the sense of distributions) of
\[
-\Div\left( \Vert \nabla g\Vert^{p-2} \nabla g
\right) = \delta_{y},\quad y\in M.
\]
A third equivalent criterion for the
$p$-hyperbolicity of $M$ is the existence of a
non-constant positive $p$-supersolution of the
$p$-Laplace equation; see \cite{H1} and
\cite{H2}. We refer to \cite{CHS}, \cite{H3}, and
\cite{HLq} for further studies on
$p$-hyperbolicity and various Liouville-type
results and to \cite{MP3} for a study of the
geometric conditions which have been previously
applied to extend the intrinsic analysis of
hyperbolicity to the extrinsic analysis which is
the main concern of the present paper.

\subsection{Outline of the paper} \label{subsecOutline}
In Section \ref{secpLap} we describe some of the
basic properties of the $p$-Laplacian and present
the corresponding maximum principle, which will
be fundamental for the comparison technique
applied in this paper. Section
\ref{secCompConstel} is devoted to set up a
so-called comparison constellation, which is
essentially molded from curvature restrictions
and a model space construction. In Section
\ref{secMain} we formulate our main result
together with three of its corollaries. They are
all proved in Sections \ref{secLapComp},
\ref{firstProofmain}, and \ref{secProofCor}. A
technical tool, the drifted 2-capacity of model
spaces is defined and analyzed in Section
\ref{secDrift}. Finally, in Section
\ref{secCapBounds} we present an alternative
proof of the main theorem based directly on
finite capacity comparison results.

\section{The $p$-Laplacian}\label{secpLap}
Let $M$ be a non-compact Riemannian manifold, with the Riemannian metric
$\langle\cdot,\cdot\rangle$ and the Riemannian volume form $d\mu$.
We say that a vector field
$\nabla u\in L^{1}_{\loc}(M)$ is a {\em{distributional gradient}} of a
function $u\in L^{1}_{\loc}(M)$ if
\[
\int_{M}\langle\nabla u,V\rangle\, d\mu  = -\int_{M}u\Div V\,d\mu
\]
for all compactly supported vector fields $V\in C^{1}_{0}(M)$. Let
$W^{1,p}(M),\ 1\le p<\infty$, be the Sobolev space of all functions
$u\in L^{p}(M)$ whose distributional gradient $\nabla u$ belongs to
$L^{p}(M)$. We equip $W^{1,p}(M)$ with the norm
$\Vert u\Vert_{1,p}=\Vert u\Vert_{p}+\Vert\nabla u\Vert_{p}$.
The corresponding local space $W^{1,p}_{\loc}(M)$ is defined in an obvious
manner. The space $W^{1,p}_{0}(M)$ is the closure of $C^{\infty}_{0}(M)$ in
$W^{1,p}(M)$.

Let $1<p<\infty$. A function $u\in W^{1,p}_{\loc}(M)$ is a (weak) solution to the
$p$-Laplace equation
\begin{equation} \label{eqPlap}
-\Div\left( \Vert \nabla u\Vert^{p-2} \nabla u
\right)  =   0
\end{equation}
in $M$ if
\begin{equation} \label{eqPlapInt}
\int_{M} \langle\Vert \nabla u\Vert^{p-2} \nabla u, \nabla
\phi \rangle \, d\mu  =   0
\end{equation}
for all $\phi \in C^{\infty}_{0}(M)$. If,
moreover, $\Vert \nabla u\Vert\in L^{p}(M)$, it
is equivalent to require \eqref{eqPlapInt} for
all $\phi\in W^{1, p}_{0}(M)$. Continuous
solutions of \eqref{eqPlap} are called
{\em{$p$-harmonic}}. Here the continuity
assumption makes no restriction since every
solution of \eqref{eqPlap} has a continuous
representative by the fundamental work of Serrin
\cite{Ser}. In fact, $p$-harmonic functions have
locally H{\"o}lder-continuous first order
derivatives by regularity results due to
Ural'tseva \cite{ural} and Lewis \cite{lew}; see
also DiBenedetto \cite{dibe}, Evans \cite{ev},
Tolksdorf \cite{tolk}, and Uhlenbeck \cite{uhl}.
Furthermore, if $D\subset M$ is a precompact open
set with $C^{1,\alpha}$ boundary ($\alpha\le 1$), 
$h\in C^{1,\alpha}(\partial D)$, and $u$ is $p$-harmonic in $D$ with 
boundary values $h$, then $u\in C^{1,\beta}(\bar D)$,
with $\beta=\beta(\alpha,p,\dim M)$, by
Lieberman \cite{lieb}.

A function $u \in W^{1, p}_{\loc}(M)$ is called a
{\em{$p$-supersolution}} in $M$ if
\begin{equation*}
\int_{M} \langle \Vert \nabla u\Vert^{p-2} \nabla u, \nabla
\phi \rangle \, d\mu  \geq  0
\end{equation*}
for all non-negative $\phi \in C^{\infty}_{0}(M)$. Similarly, a
function $v \in W^{1, p}_{\loc}(M)$ is called a
{\em{$p$-subsolution}} in $M$ if
\begin{equation*}
\int_{M} \langle \Vert \nabla v\Vert^{p-2}\nabla v, \nabla
\phi \rangle \, d\mu  \leq  0
\end{equation*}
for all non-negative $\phi \in C^{\infty}_{0}(M)$. A fundamental feature of
solutions of \eqref{eqPlap} is the following {\em{maximum (or comparison)
principle}} which will be instrumental for the comparison technique
presented below in Sections \ref{secMain} and \ref{secDrift}:
If $u\in W^{1,p}(M)$ is a $p$-supersolution,
$v\in W^{1,p}(M)$ is a $p$-subsolution, and $\max(v-u,0)\in W^{1,p}_{0}(M)$,
then $u\ge v$ a.e. in $M$. In particular, if $D\subset M$ is a
precompact open set, $u\in C(\bar D)$ is a $p$-supersolution,
$v\in C(\bar D)$ is a $p$-subsolution, and $u\ge v$ in $\partial D$,
then $u\ge v$ in $D$. For the reader's convenience we recall the short proof
of the comparison principle from \cite[3.18]{HKM}. The proof is based on the 
following elementary inequality:
Let $a \neq b$ denote two vectors in a given tangent space $T_{x}M$
and suppose that $1 < p < \infty$. Then
\[
\langle \Vert a \Vert^{p-2}a - \Vert b \Vert^{p-2}b,  a - b
\rangle >  0.
\]
Suppose then that $u\in W^{1,p}(M)$ is a $p$-supersolution and
$v\in W^{1,p}(M)$ is a $p$-subsolution such that 
$\phi=\max(v-u,0)\in W^{1,p}_{0}(M)$.
Since
\begin{align*}
0&\ge \int_{M}\langle\Vert\nabla v\Vert^{p-2}\nabla
v,\nabla\phi\rangle\,d\mu
-\int_{M}\langle\Vert\nabla u\Vert^{p-2}\nabla
u,\nabla\phi\rangle\,d\mu\\
&=\int_{\{u<v\}}\langle\Vert\nabla v\Vert^{p-2}\nabla v
-\Vert\nabla u\Vert^{p-2}\nabla u,\nabla v-\nabla u\rangle\,d\mu\ge 0,
\end{align*}
we have $\nabla\phi =0$ a.e. in $M$ as required.

\section{Comparison Constellations} \label{secCompConstel} 
We assume
throughout the paper that $S^{m}$ is a non-compact
connected complete Riemannian submanifold of a complete
Riemannian manifold $N^{n}$. Furthermore, we
assume that $N^{n}$
possesses at least one pole. Recall that a pole
is a point $o$ such that the exponential map
$\exp_{o}\colon T_{o}N^{n} \to N^{n}$ is a
diffeomorphism. For example, an Hadamard--Cartan
manifold has everywhere non-positive sectional
curvatures and since it is also by definition
simply connected, {\em{every}} point is a pole.
The r\^{o}le of the pole $o$ is precisely to
serve as the origin of a smooth distance function
$r$ from $o$: For every $x \in N^{n}\setminus \{o\}$ we
define $r(x) = \dist_{N}(o, x)$, and this
distance is realized by the length of a unique
geodesic from $o$ to $x$, which is the {\it
radial geodesic from $o$}. We also denote by $r$
the restriction $r\vert_S: S\to \erre_{+} \cup
\{0\}$. This restriction is called the
{\em{extrinsic distance function}} from $o$ in
$S^m$. The gradients of $r$ in $N$ and $S$ are
denoted by $\nabla^N r$ and $\nabla^S r$,
respectively. Let us remark that $\nabla^S r(x)$
is just the tangential component in $S$ of
$\nabla^N r(x)$, for all $x\in S$. Then we have
the following basic relation:
\begin{equation*}
\nabla^N r = \nabla^S r +(\nabla^N r)^\bot ,
\end{equation*}
where $(\nabla^N r)^\bot(x)$ is perpendicular to
$T_{x}S$ for all $x\in S$.

\subsection{Curvature restrictions}
\label{subsecCurvRestrict} The sectional
curvatures of $N$ {\em{along}} the radial
geodesics from $o$ are called the $o$-radial
sectional curvatures of $N$.

\begin{definition}\label{defRadCurv}
Let $o$ be a point in a Riemannian manifold $M$
and let $x \in M\setminus\{ o \}$. The sectional
curvature $K_{M}(\sigma_{x})$ of the two-plane
$\sigma_{x} \in T_{x}M$ is then called an
\textit{$o$-radial sectional curvature} of $M$ at
$x$ if $\sigma_{x}$ contains the tangent vector
to a minimal geodesic from $o$ to $x$. We denote
these curvatures by $K_{o, M}(\sigma_{x})$.
\end{definition}

The $o$-radial sectional curvatures of $N$
control the second order behavior of $r(x)$ in
$N$ via the classical Jacobi field index theory.
Indeed, a bound on the $o$-radial sectional
curvatures gives a bound on the Hessian of radial
functions, $\Hess^{N}\left(f(r)\right)$, as
proved by Greene and Wu \cite[Theorem A]{GreW};
see Theorem \ref{thmGreW} below. The submanifold
$S$ and the restricted radial functions
$f(r)|_{S}$ inherit this second order bound to
the $S$-intrinsic Hessian, $\Hess^{S}f(r)$, and
therefore also to the Laplacian $\Delta^{S}f(r)$
of such modified distance functions.

The mean curvatures $H_{S}$ of $S$ also appear in
the Laplacian $\Delta^{S}f(r)$ via its radially
weighted component, which we define as follows:

\begin{definition}
The {\em{$o$-radial mean convexity}}
$\,\mathcal{C}(x)$ of $S$ in $N$, is defined in
terms of the inner product of $H_{S}$ with the
$N$-gradient of the distance function $r(x)$ as
follows:
\begin{equation*}
\mathcal{C}(x) \,=\, -\langle \nabla^{N}r(x),
H_{S}(x) \rangle, \quad x \in S,
\end{equation*}
where $H_{S}(x)$ denotes the mean curvature
vector of $S$ in $N$, i.e. the mean trace of the
second fundamental form $\alpha_{x}$. With
respect to an orthonormal basis $\{X_{1}, ...,
X_{m}\}$ of $T_{x}S$ at $x \in S$ we have
\begin{equation*}
H_{S}(x) =  \frac{1}{m}\sum_{i=1}^{m}
\alpha_{x}\left(X_{i}, X_{i}\right).
\end{equation*}
\end{definition}

We will assume, that $\mathcal{C}(x)$ is bounded
from above by a function $h(r(x))$ which only
depends on the distance $r$ from $o\, $:
\begin{equation*}
\mathcal{C}(x) \leq h(r(x)), \quad x \in S.
\end{equation*}

Moreover, for $p > 2$ we shall also need a
particular inequality for the second fundamental
form of $S$ in $N$ in the direction of the
gradient $\nabla^{N}r(x)$. This gives rise to
the following definition:

\begin{definition}\label{defBetax}
The {\em{$o$-radial component $\,\mathcal{B}(x)$
of the second fundamental form}}  of $S$ in $N$,
is defined in terms of the following inner
product:
\begin{equation*}
\mathcal{B}(x)  =  - \langle\nabla^{N}r(x),
\alpha_{x}(U_{r}, U_{r})\rangle,
\end{equation*}
where
\begin{equation*}
U_{r} =  \nabla^{S}(r(x))/\Vert
\nabla^{S}r(x) \Vert \in T_{x}S \subset T_{x}N
\end{equation*}
is the unit tangent vector to $S$ in the
direction of $\nabla^{S}r(x)$ (resp. tacitly
assumed to be $0$ in case $\nabla^{S}r(x) = 0$).
\end{definition}

We assume that $\mathcal{B}(x) $ is bounded from
above by a function $\lambda(r(x))$ which only
depends on the distance $r$ from $o\, $:
\begin{equation*}
\mathcal{B}(x) \leq  \lambda(r(x)).
\end{equation*}

Finally, we also impose an upper control on the
'radiality' of the submanifold, i.e. a local
measure of how much the submanifold is extending
away from the pole $o$:

\begin{definition}
The {\em{\,$o$-radial tangency}} $\mathcal{T}(x)$
of $S$ in $N$ is defined as follows:
\begin{equation*}
\mathcal{T}(x)  =  \Vert \nabla^{S}r(x)\Vert
\end{equation*}
for all $x\in S$.
\end{definition}

We assume that this $S$-gradient of the
restricted distance function $r|_{S}$ has an
upper radial support function $g(r) \leq 1$:
\begin{equation*}
\mathcal{T}(x)   \leq  g(r(x)) .
\end{equation*}

\begin{definition}
Given a connected and complete
$m$-dimen\-sional submanifold $S^m$ in a complete
Riemannian manifold $N^n$ with a pole $o$, we
denote
the {\em{extrinsic metric balls}} of
(sufficiently large) radius $R$ and center $o$ by
$D_R(o)$. They are defined as any connected
component of the intersection
$$
B_{R}(o) \cap S =\{x\in S \colon r(x)< R\},
$$
where $B_{R}(o)$ denotes the open geodesic ball
of radius $R$ centered at the pole $o$ in
$N^{n}$.
Using these extrinsic balls we define the
$o$-centered extrinsic annuli
$$
A_{\rho,R}(o)= D_R(o) \setminus \bar D_{\rho}(o)
$$
in $S^m$ for $\rho < R$, where $D_{R}(o)$ is the component
of $B_{R}(o) \cap S$ containing $D_{\rho}(o)$.
\end{definition}

The upper bounding functions
$h(r)$, $g(r)$, and $\lambda(r)$ together with a
suitable control on the $o$-radial sectional
curvatures of the ambient space will eventually
control the $p$-Laplacian of restricted radial
functions on $S$. In particular, we consider
potential functions stemming from capacity
calculations of radially symmetric comparison
spaces and transplant them to $S$ via the
distance function $r$ in $N$. Such
transplantations are then compared with the
'correct' potentials on extrinsic metric balls of
$S$. The maximum principle for the $p$-Laplacian
$\pL^{S}$ then finally gives the comparison
result for capacities in $S$. Concerning the
general strategy and types of results (in the
case of $p =2$) we refer to \cite{M1},
\cite{Palmer1}, and \cite{MP1}.

We now collect the previous ingredients and
formulate the general framework for our
$p$-hyperbolicity comparison result:

\begin{definition}\label{defConstellat}
Let $N^{n}$ denote a Riemannian manifold with a
pole $o$ and distance function $r = r(x) =
\dist_{N}(o, x)$. Let $S^{m}$ denote a connected
complete submanifold in $N^{n}$ and assume that
there is an extrinsic ball $D_{\rho}(o)$ 
which is precompact with
smooth boundary $\partial D_{\rho}(o)$ in
$S^{m}$.
 Let $M_{w}^{m}$
denote a $w$-model with warping function $w:
\pi(M_{w}^{m}) \to \mathbb{R}_{+}$ and center
$o_{w}$; see Definition \ref{defModel}. Then the triple $\{ N^{n}, S^{m},
M_{w}^{m} \}$ is called a {\em{comparison
constellation}} on the interval $[0, R]$ if the
$o$-radial sectional curvatures of $N$ are
bounded from above by the $o_{w}$-radial
sectional curvatures of $M_{w}^{m}$:
\begin{equation} \label{eqKcomp}
K_{o, N}(\sigma_{x}) \leq
-\frac{w''(r)}{w(r)}
\end{equation}
for all $x$ with $r=r(x)\in [0,R]$ and, moreover,
the radial tangency $\mathcal{T}$ and the radial
convexity functions  $\mathcal{B}$ and
$\mathcal{C}$ of the submanifold $S^{m}$ are all
bounded from above by smooth radial functions
$g(r)$, $\lambda(r)$, and $h(r)$, respectively:
\begin{equation}
\label{eqRadConv}
\begin{aligned}
\mathcal{T}(x) &\leq g(r(x)),\\
\mathcal{B}(x) &\leq \lambda(r(x)) ,\,\,\, \textrm{and}\\
\mathcal{C}(x) &\leq h(r(x)) \,\,\,\textrm{for
all}\,\,\, x \in S^{m} \,\,\,
{\textrm{with}}\,\,\, r(x) \in [0, R] .
 \end{aligned}
\end{equation}
\end{definition}

\begin{remark}\label{theRemk}
This definition of comparison constellation
extends a previous definition considered in
\cite{MP3}. In that paper, the triple $\{ N^{n},
S^{m}, M_{w}^{m} \}$ is called a {\em{comparison
constellation}} if inequality (\ref{eqKcomp}) holds
and if in addition only the following condition
holds in replacement of inequalities
(\ref{eqRadConv}) for some bounding radial
function $h(r)$:
\begin{equation*}
\mathcal{C}(x) \leq h(r(x)) \leq
\frac{w'(r(x))}{w(r(x))} \,\,\,\textrm{for
all}\,\,\, x \in S^{m}.
\end{equation*}
It is proved in \cite{MP3} that under these
conditions $S^m$ is $2$-hyperbolic if
\begin{equation*}
\int_{\rho}^{\infty}\frac{\mathcal{G}^m(r)}{w^{m-1}(r)}dr
<\infty ,
\end{equation*}
where
\begin{equation*}
\mathcal{G}(r) \,=\, \exp(\int_{\rho}^r h(t) dt).
\end{equation*}
\end{remark}

\subsection{Warped products and model spaces}\label{subsecWarp}
Warped products are generalized manifolds of revolution, see e.g.
\cite{O'N}. Let $(B^{k}, g_{B})$ and $(F^{l}, g_{F})$ denote two
Riemannian manifolds and let $w \colon B \to \mathbb{R_{+}}$ be a
positive real function on $B$. We assume throughout that $w$ is
at least $C^{1}$ with piecewise continuous second order derivatives.
We consider the product manifold $M^{k+l}= B \times F $
and denote the projections onto the factors by $\pi\colon
M \to B$ and $\sigma\colon M \to F$, respectively. The metric
$g$ on $M$ is then defined by the following $w$-modified
(warped) product metric
\begin{equation*}
 g = \pi^{*}(g_{B}) + (w \circ \pi)^{2} \sigma^{*}(g_{F}).
\end{equation*}
\begin{definition}
The Riemannian manifold $(M, g) = (B^{k} \times F^{l}, g) $ is
called a \textit{warped product} with \textit{warping function} $w$,
base manifold $B$ and fiber $F$. We write as follows: $M_{w}^{m} =
B^{k} \times_{w} F^{l}$.
\end{definition}

\begin{definition}[See \cite{Gri}, \cite{GreW}]\label{defModel}
A $w-$model $M_{w}^{m}$ is a smooth warped product with base $B^{1}
= [0,\Lambda[ \,\subset \mathbb{R}$ (where $0 < \Lambda
\leq  \infty$), fiber $F^{m-1} = \mathbb{S}^{m-1}_{1}$ (i.e. the unit
$(m-1)$-sphere with standard metric), and warping function $w\colon
[0,\Lambda[ \to \mathbb{R}_{+}\cup \{0\}$, with $w(0) = 0$,
$w'(0) = 1$, and $w(r) > 0$ for all $r >  0$. The point
$o_{w} = \pi^{-1}(0)$, where $\pi$ denotes the projection onto
$B^1$, is called the {\em{center point}} of the model space. If
$\Lambda = \infty$, then $o_{w}$ is a pole of $M_{w}^{m}$.
\end{definition}

\begin{proposition}\label{propSpaceForm}
The simply connected space forms $\mathbb{K}^{m}(b)$ of constant
curvature $b$ are $w-$models with warping functions
\begin{equation*}
w(r) = Q_{b}(r) =\begin{cases} \frac{1}{\sqrt{b}}\sin(\sqrt{b}\, r) &\text{if $b>0$}\\
\phantom{\frac{1}{\sqrt{b}}} r &\text{if $b=0$}\\
\frac{1}{\sqrt{-b}}\sinh(\sqrt{-b}\,r) &\text{if $b<0$}.
\end{cases}
\end{equation*}
Note that for $b > 0$ the function $Q_{b}(r)$ admits a smooth
extension to  $r = \pi/\sqrt{b}$.
\end{proposition}

\begin{proposition}[See e.g. \cite{O'N}]\label{propModelDist}
Let $M_{w}^{m} = B^{1} \times_{w} \mathbb{S}_{1}^{m-1}$ be a $w-$model. Let
$r_{0}$ and $r$ denote two points in $B^{1}$. Then the geodesic distance
from every $x \in \pi^{-1}(r)$ to $\pi^{-1}(r_{0})$ is
$|r - r_{0}|$.
\end{proposition}

\begin{proposition}[See \cite{O'N} p. 206]\label{propWarpMean}
Let $M_{w}^{m}$ be a $w-$model with warping function $w(r)$ and
center $o_{w}$. The distance sphere of radius $r$ and center $o_{w}$
in $M_{w}^{m}$ is the fiber $\pi^{-1}(r)$. This distance sphere has
the following constant mean curvature vector in $M_{w}^{m}$
\begin{equation*}
H_{{\pi^{-1}}(r)} = -\eta_{w}(r)\,\nabla^{M}\pi =
-\eta_{w}(r)\,\nabla^{M}r,
\end{equation*}
where the mean curvature function $\eta_{w}(r)$ is defined by
\begin{equation*}
\eta_{w}(r)  = \frac{w'(r)}{w(r)} = \frac{d}{dr}\log (w(r)).
\end{equation*}
\end{proposition}

In particular, we have for the constant curvature space forms
$\mathbb{K}^{m}(b)$:
\begin{equation*}
\eta_{Q_{b}}(r) = \begin{cases} \sqrt{b}\cot(\sqrt{b}\,r) &\text{if $b>0$}\\
\phantom{\sqrt{b}} 1/r &\text{if
$b=0$}\\\sqrt{-b}\coth(\sqrt{-b}\,r) &\text{if $b<0$}\,.
\end{cases}
\end{equation*}

The radial curvature in model spaces is given by the following result

\begin{proposition}[See \cite{GreW} and \cite{Gri}] \label{propModelRadialCurv}
Let $M_{w}^{m}$ be a $w-$model with center point $o_{w}$. Then the
$o_{w}$-radial sectional curvatures of $M_{w}^{m}$ at every $x \in
\pi^{-1}(r)$ (for $r > 0$) are all identical and determined
by
\begin{equation*}
K_{o_{w} , M_{w}}(\sigma_{x}) = -\frac{w''(r)}{w(r)}.
\end{equation*}
\end{proposition}

\subsection{Hessian and Laplacian comparison analysis}\label{subsecLap}
Concerning the second order analysis of the distance function
$r$ we need firstly and foremost the Hessian comparison theorem
for manifolds with a pole:

\begin{theorem}[See \cite{GreW}, Theorem A]\label{thmGreW}
Let $N=N^{n}$ be a manifold with a pole $o$, let $M=M_{w}^{m}$ denote a
$w-$model with center $o_{w}$, and $m \leq n$. Suppose that every $o$-radial
sectional curvature at $x \in N \setminus \{o\}$ is bounded from above by
the $o_{w}$-radial sectional curvatures in $M_{w}^{m}$ as follows:
\begin{equation*}
K_{o, N}(\sigma_{x}) \leq -\frac{w''(r)}{w(r)}
\end{equation*}
for every radial two-plane $\sigma_{x} \in T_{x}N$ at distance $r =
r(x) = \dist_{N}(o, x)$ from $o$ in $N$. Then the Hessian of the
distance function in $N$ satisfies
\begin{equation}\label{eqHess}
\begin{aligned}
\Hess^{N}(r(x))(X, X) &\geq \Hess^{M}(r(y))(Y, Y)\\ &=
\eta_{w}(r)\left(1 - \langle \nabla^{M}r(y), Y \rangle_{M}^{2}
\right) \\ &= \eta_{w}(r)\left(1 - \langle \nabla^{N}r(x), X
\rangle_{N}^{2} \right)
\end{aligned}
\end{equation}
for every unit vector $X$ in $T_{x}N$ and for every unit vector $Y$
in $T_{y}M$ with $\,r(y) = r(x) = r\,$ and $\, \langle
\nabla^{M}r(y), Y \rangle_{M} = \langle \nabla^{N}r(x), X
\rangle_{N}\,$.
\end{theorem}

\begin{remark}
In \cite[Theorem A, p. 19]{GreW}, the Hessian of
$r_M$ is less or equal to the Hessian of $r_N$
provided that the radial curvatures of $N$ are
bounded from above by the radial curvatures of
$M$ and provided that $\dim M \geq \dim N$. This
latter dimension condition is {\em{not}}
satisfied in our setting. However, since $(M^{m},
g)$ is a $w-$model space it has an
$n-$dimensional $w-$model space companion with
the same radial curvatures and the same Hessian
of radial functions as $(M^{m}, g)$. In effect,
therefore, applying \cite[Theorem A, p. 19]{GreW}
to the high-dimensional comparison space gives
the low-dimensional comparison inequality as
stated.
\end{remark}

If $\mu\colon N \to \mathbb{R}$ denotes a smooth
function on the ambient space $N$, then the restriction $\Tilde{\mu}= \mu_{|_{S}}$
is a smooth function on the submanifold $S$ and the respective Hessian tensors,
$\operatorname{Hess}^{N}(\mu)$ and
$\operatorname{Hess}^{S}(\Tilde{\mu})$,
are related as follows:

\begin{proposition}[\cite{JK}]
\begin{equation} \label{eqHessOrig}
 \operatorname{Hess}^{S}(\Tilde{\mu})(X, Y)=
    \operatorname{Hess}^{N}(\mu)(X, Y)
     + \langle \nabla^{N}(\mu),
\alpha_{x}(X,Y)\rangle
\end{equation}
for all tangent vectors $X , \, Y  \in T_{x}S^{m} \subset
T_{x}N^{n}$, where $\alpha_{x}$ is the second fundamental
form of $S$ at $x$ in $N$.
\end{proposition}

If we compose $\mu$ with a smooth function $f\colon \mathbb{R} \to
\mathbb{R}$ we then get:
\begin{corollary}[\cite{JK}] \label{corHessfomu}
\begin{equation*}
 \begin{aligned}
  \operatorname{Hess}^{S}(f \circ \Tilde{\mu})(X, X)
&  =  f''(\mu) \langle\nabla^{N}(\mu), X \, \rangle^{2}\\
 + f'(\mu)\bigl( \operatorname{Hess}^{N}(\mu)&(X, X)
 +  \langle \nabla^{N}(\mu),  \alpha_{x}(X,  X) \rangle
\bigr)
 \end{aligned}
 \end{equation*}
 for all $X \in T_{x}S^{m}$.
\end{corollary}

Combining the estimate \eqref{eqHess} with Corollary \ref{corHessfomu} and tracing
the resulting Hessian comparison
statement in an orthonormal basis of $T_{x}S^{m}$, we obtain the
following instrumental inequality for the Laplacian of
(extrinsic) radial functions restricted to the submanifold $S$:
\begin{proposition} \label{corLapComp}
Suppose that the assumptions of Theorem
\ref{thmGreW} are satisfied. Then we have for
every smooth real-valued function $f\circ r$ with $f' \geq 0$
the following inequality for the standard Laplacian:
\[
\Delta^{S}(f \circ r) \geq  \left(f''(r)
- f'(r)\eta_{w}(r) \right)
 \Vert \nabla^{S} r \Vert^{2}
+ m f'(r) \left(\eta_{w}(r) +\langle\nabla^{N}r,H_{S}\rangle\right),
\]
where $H_{S}$ denoted the mean curvature vector of $S$ in $N$.
\end{proposition}

\section{Main results}\label{secMain}
Applying the notion of a comparison constellation
as defined in the previous section, we now
formulate our main $p$-hyperbolicity result. The
proofs are developed through the following
sections.

\begin{theorem} \label{thmMain}
Consider a comparison constellation
$\{N^{n},S^{m},M_{w}^{m}\}$ on the interval
$[\,0, \infty[\,$. Assume further that the
functions $h(r)$ and $\lambda(r)$ are
{\em{balanced}} with respect to the warping
function $w(r)$ by the following inequality:
\begin{equation} \label{eqBalance}
\mathcal{M}(r) :=  \left( m + p -
2\right)\eta_{w}(r) - m\,h(r) - (p-2)\lambda(r)
\geq  0.
\end{equation}
Let $\Lambda(r)$ denote the function
\begin{equation*}
\Lambda(r) =
w(r)\exp\left(-\int_{1}^{r}
\frac{\mathcal{M}(t)}{(p-1)g^{2}(t)}\,
dt\right).
\end{equation*}
Suppose finally that $p\ge 2$ and that
\begin{equation} \label{eqTransienceCond}
\lim_{R \to \infty}  \int_{\rho}^{R} \Lambda(t)
\,  dt  <  \infty.
\end{equation}
Then $S^{m}$ is $p$-hyperbolic.
\end{theorem}

We observe the following corollaries; they will
be proved in Section \ref{secProofCor}.

\begin{corollary} \label{corHadamard}
Suppose (in Theorem \ref{thmMain}) that we can
choose $w(r)=  Q_{b}(r) =
\sinh(\sqrt{-b}\,r)/\sqrt{-b}$ for some $b<0$,
i.e. we apply the negatively curved space form
$\mathbb{K}^{m}(b)$ to play the role of a model
space in the comparison constellation. Suppose
that there exist constants $\lambda_{0}$ and
$h_{0}$ such that
\begin{equation*}
\begin{aligned}
\mathcal{B}(x) &\leq \lambda_{0}\, \, \, \textrm{and}\\
\mathcal{C}(x) &\leq h_{0} \,\,\,\textrm{for
all}\,\,\, x \in S^{m}.
 \end{aligned}
\end{equation*}
Suppose further that for some $\tilde p \geq 2$ we have
\begin{equation}\label{eqCorHad}
m\,h_{0} + (\tilde p-2)\lambda_{0} <
(m-1)\sqrt{-b}.
\end{equation}
Then $S^{m}$ is $p$-hyperbolic for all $p$ in the
range $2\leq p\leq \tilde p$.
\end{corollary}

\begin{corollary}\label{corTroy}
Consider a purely intrinsic setting and comparison
constellation: $S^{n} = N^{n} = M_{w}^{n}$.
 Then $S^{n}$ is
$p$-hyperbolic if and only if
\begin{equation*}
\label{eqTroy}
\int_{\rho}^{\infty}\,\frac{1}{w(t)^{\frac{n-1}{p-1}}}\,dr
<\infty.
\end{equation*}
\end{corollary}
This observation is originally due to M.
Troyanov, see \cite[Corollary 5.4]{T1}.

\begin{corollary} \label{exTroy}
Let $(M^{m}, g)$ denote a complete  manifold with
intrinsic concentric metric balls $B_{r}(o)$
centered at $o \in M$. Suppose that for some
$p\geq 2$ and for some $\rho >0$ we have
\begin{equation} \label{eqTroyInfinit}
\int_{\rho}^{\infty}\,\frac{1}{\Vol(\partial
B_{r}(o))^{\frac{1}{p-1}}}\,\,dr =  \infty,
\end{equation}
and suppose that there are constants $\lambda_{0}
> 0 $ and $b < 0$ so that
\begin{equation} \label{eqplm}
(p-2)\lambda_{0} \, < (m - 1)\,\sqrt{-b}.
\end{equation}
Then $(M, g)$ does not admit a minimal isometric
immersion with bounded second fundamental form
$\Vert \alpha\Vert \, \leq \, \lambda_{0}$ into
any Hadamard--Cartan  manifold $N^{n}$, $n \geq
m$,  with sectional curvatures bounded from above
by $b$.

\begin{proof}
Condition (\ref{eqTroyInfinit}) implies that the
manifold $(M^{m}, g)$ is $p$-parabolic according to
\cite[Corollary 5.4]{T1}, whereas the condition
(\ref{eqplm}) implies $p$-hyperbolicity of $(M^{m},
g)$ according to Corollary \ref{corHadamard} of
the present work - upon observing that 
$\mathcal{C}(x)\equiv 0$ by the minimality assumption 
and that 
$\,\Vert \alpha_{x}\Vert \, \leq \, \lambda_{0}\,$ 
implies
$\,\mathcal{B}(x)\,\leq \, \lambda_{0}\,$.

\end{proof}
\end{corollary}

\section{Drifted $2$-capacity of model spaces} \label{secDrift}

\begin{definition} \label{defLapmod}
Let $(M, g)$ denote a Riemannian manifold with
Laplace operator $\Delta^{M}$, and let $V$ denote
a continuous vector field on $M$. The drifted
Brownian motion on $M$ with the drift vector field
$V$ is then generated by the modified Laplacian
$\Lmod$
\begin{equation*}
\Lmod f =  \Delta^{M}f + \langle  \nabla^{M} f, V \rangle
\end{equation*}
for every smooth function $f$ on $M$.
\end{definition}

We consider, in particular, the drift vector field
\begin{equation*}
V = \mathcal{V}(r)\nabla^{M}r
\end{equation*}
with
\begin{equation*}
\mathcal{V}(r) =
\frac{\mathcal{M}(r)}{(p-1)g^{2}(r)} -
m\,\eta_{w}(r)
\end{equation*}
on model spaces $M = M^{m}_{w}$, so that the
modified Laplacian then reads as
\begin{equation*}
\Lmod \psi(x) =  \Delta^{M}\psi(x) +
\psi'(r(x))\mathcal{V}(r(x))
\end{equation*}
for smooth functions $\psi$ on $M^{m}_{w}$. For
purely {\em{radial}} functions $\psi(r)$ we get
\begin{lemma}\label{lemLrad}
Let $\psi = \psi(r)$ denote a function on the
$w$-model space $M =M_{w}^{m}$ which only depends
on the radial distance $r$ to the center $o_{w}$.
Then
\begin{equation*}
\Lmod \psi(r) = \psi''(r) + \psi'(r)\left(
\frac{\mathcal{M}(r)}{(p-1)\,g^{2}(r)} -
\eta_{w}(r)\right).
\end{equation*}
\end{lemma}

The Dirichlet problem associated to $\Lmod$
defined on so-called {\em{extrinsic annuli}} is
defined as follows:

 First, the annular domains in the model
space are denoted by
\[
A_{\rho, R}^{w} = \{x \in M^{n}_{w}\colon
\pi(x) \in [\rho, R] \} = \pi^{-1}([\rho, R]),
\]
and the corresponding boundaries are denoted by
$\partial D_{\rho}^{w} = \pi^{-1}(\rho)$ and
$\partial D_{R}^{w} = \pi^{-1}(R)$, respectively.
We consider the unique radial function $\psi_{\rho,R}(r)$
which solves the one-dimensional
Laplace-Dirichlet problem on the model space
annulus $A_{\rho, R}^{w}$:
\begin{equation}\label{eqDir}
\begin{cases}
\Lmod \psi &= 0\,\,\,\text{on $A_{\rho, R}^{w}$}\\
\phantom{L }\psi &=0\,\,\,\text{on $\partial
D_{\rho}^{w}$}  \\
\phantom{L }\psi&= 1\,\,\,\text{on $\partial
D_{R}^{w}$}.
\end{cases}
\end{equation}

The explicit solution to the Dirichlet problem
\eqref{eqDir} is given in the following
Proposition,  with a focus towards the
corresponding expression for the drifted annular
capacity in the model space; see \cite{MP3},
\cite{MP2}, and Section \ref{secCapBounds} below.

\begin{proposition}\label{propDirSol}
The solution to the Dirichlet problem (\ref{eqDir}) only depends on
$r$ and is given explicitly - via the function
$\Lambda(r)$ introduced in Theorem \ref{thmMain}, by:
\begin{equation} \label{eqPsi}
\psi_{\rho,R}(r) =  \frac{\int_{\rho}^{r}
\Lambda(t)\,dt}{\int_{\rho}^{R}\Lambda(t)\,dt}.
\end{equation}
The corresponding 'drifted' $2$-capacity is
\begin{equation} \label{eqModelCap}
\begin{aligned}
\C_{\Lmod}(A_{\rho, R}^{w})
&=\int_{\partial D_{\rho}^{w}}\langle\nabla^{M}\psi_{\rho,R},\nu\rangle\,dA\\
&=\Vol(\partial D_{\rho}^{w})\Lambda(\rho)\left(\int_{\rho}^{R}
\Lambda(t)\,dt\right)^{-1}.
\end{aligned}
\end{equation}
\end{proposition}

\section{$p$-Laplacian comparison}\label{secLapComp}

Let us consider comparison constellations
$\{N^{n},S^{m},M_{w}^{m}\}$ on intervals $[0,R]$
for $R>0$. Since the  $o$-radial mean convexity
of $S$ has an upper bound
\[
\mathcal{C}(x) = -\langle \nabla^{N}r(x),
H_{S}(x) \rangle  \leq h(r(x)),
\]
we obtain the following estimate using Proposition  \ref{corLapComp}
\begin{equation}
\label{eqLap1}
\Delta^{S}(f \circ r) \geq  \left(f''(r)
- f'(r)\eta_{w}(r) \right)
 \Vert \nabla^{S} r \Vert^{2}
+ m\,f'(r) \left(\eta_{w}(r) - h(r)\right).
\end{equation}

In what follows we use shorthand
$F(x)=f'(r(x))\Vert \nabla^{S} r(x)\Vert\,$ for
all $x\in S$ to simplify the notation. To get
estimates for the $p$-Laplacian of $f\circ r$ we
first compute
\begin{align*}
&\pL^{S}f(r(x)) = \Div^{S}\left( \Vert\nabla^{S}
f(r(x))\Vert^{p-2}\,\nabla^{S} f(r(x))
\right) \\
&= \Vert \nabla^{S} f(r(x))\Vert^{p-2}\,\Delta^{S}f(r(x))
+ \left\langle \nabla^{S}\Vert
\nabla^{S} f(r(x))\Vert^{p-2} ,  \nabla^{S}f(r(x))
\right\rangle \\
&=F^{p-2}(x)\Delta^{S}f(r(x))
+ \left\langle
\nabla^{S}F^{p-2}(x) , f'(r(x))\nabla^{S}r(x)\right\rangle
\\
&= F^{p-2}(x)\Delta^{S}f(r(x))\\
&\phantom{mm}+ \Bigl\langle(p-2)F^{p-3}(x)\left(f''(r(x))\Vert
\nabla^{S}r(x) \Vert\nabla^{S}r(x) +
f'(r(x))\nabla^{S}\Vert \nabla^{S}r(x)
\Vert \right),\\
&\phantom{mmmm} f'(r(x))\nabla^{S}r (x)\Bigr\rangle \\
&= F^{p-2}(x) \Bigl((p-2) \Bigl(
f''(r(x))\Vert \nabla^{S} r(x)\Vert^{2} +
f'(r(x))\frac{\left\langle \nabla^{S}r(x), \nabla^{S}\Vert
\nabla^{S}r(x) \Vert\right \rangle}{\Vert \nabla^{S}r(x)
\Vert}\Bigr)\\
&\phantom{mmm} + \Delta^{S}f(r(x))\Bigr).
\end{align*}
This partial 'isolation' of the factor $(p-2)$ is
the reason behind the general assumption $p \geq
2$ in this work.  The factor on $(p-2)$ is
controlled via the following observation, which
introduces the bound $\lambda(r)$ into this
setting:
\begin{lemma} \label{lemLambda}
Let $\{N^{n},S^{m},M_{w}^{m}\}$ be a comparison 
constellation on $[0,R]$ for $R>0$. Suppose that the 
$o$-radial component of the second fundamental form
of $S$ (see Definition \ref{defBetax}) has an upper bound 
\[
\mathcal{B}(x)  \le \lambda(r(x)).
\]
Then
\begin{equation} \label{eqLambda}
\begin{aligned}
&\frac{\left \langle \nabla^{S}r(x), \nabla^{S}\Vert
\nabla^{S}r(x)\Vert\right \rangle }{\Vert
\nabla^{S}r(x) \Vert}\, \\ &\qquad \qquad \qquad
=  \Hess^{S}(r(x))\left(U_{r},
U_{r}\right) \\
& \qquad  \qquad \qquad =
\Hess^{N}(r(x))\left(U_{r}, U_{r}\right) +
\left\langle \nabla^{N}r(x),  \alpha_{x}\left(U_{r},
U_{r} \right)\right \rangle \\
&\qquad  \qquad \qquad \geq  \eta_{w}(r(x))\left(
1 - \Vert \nabla^{S}r(x)\Vert^{2} \right) -
\lambda(r(x)) .
\end{aligned}
\end{equation}
\end{lemma}
\begin{proof}
By definition of the Hessian via the induced
connection $\D^{S}$ in $S$ we have directly for
the first equality in (\ref{eqLambda}):
\begin{equation*}
\begin{aligned}
\Hess^{S}(r)\left(\nabla^{S}r, \nabla^{S}r\right)
&= \left \langle \D^{S}_{\nabla^{S}r}\nabla^{S}r,
\nabla^{S}r \right\rangle \\
&=  \tfrac{1}{2}\D^{S}_{\nabla^{S}r}\left\langle
\nabla^{S}r, \nabla^{S}r
\right\rangle \\
&= \tfrac{1}{2}{\nabla^{S}r}\left\langle \nabla^{S}r, \nabla^{S}r \right\rangle \\
&= \tfrac{1}{2}\left\langle \nabla^{S}\Vert \nabla^{S}r\Vert^{2}, \nabla^{S}r \right\rangle\\
&= \Vert \nabla^{S}r \Vert\left\langle
\nabla^{S}\Vert \nabla^{S}r\Vert, \nabla^{S}r
\right\rangle ,
\end{aligned}
\end{equation*}
so that
\begin{equation*}
\begin{aligned}
\Hess^{S}(r(x))\left(U_{r}, U_{r}\right)  &=
\frac{\Hess^{S}(r)\left( \nabla^{S}r,
\nabla^{S}r\right)}{\Vert \nabla^{S}r
\Vert^{2}} \\
&=\frac{\left \langle \nabla^{S}r(x), \nabla^{S}\Vert
\nabla^{S}r(x)\Vert\right \rangle }{\Vert
\nabla^{S}r(x) \Vert}.
\end{aligned}
\end{equation*}
The other (in)equalities in (\ref{eqLambda})
follow from (\ref{eqHessOrig}) and
(\ref{eqHess}), respectively.
\end{proof}

The following result relates the $p$-Laplacian of a
radial function $f(r)$ with its $2$-drifted
Laplacian, as defined in Section \ref{secDrift}.

\begin{lemma}\label{lemPLapComp}
Let $\{N^{n},S^{m},M_{w}^{m}\}$ be a comparison constellation on $[0,R]$
for $R>0$. Let $f\circ r$ be a smooth real-valued function with $f' \geq 0$,
and suppose now that $f(r)$ satisfies the following
condition (to be molded shortly from the balance
condition (\ref{eqBalance})):
\begin{equation}\label{parenth}
f''(r) - f'(r)\eta_{w}(r) \leq  0.
\end{equation}
Then, for all $x \in S$,
\[
\pL^{S}f(r(x)) \geq (p-1)F^{p-2}(x)g^{2}(r(x))
\Lmod(f(r(x))) ,
\]
where $L$ is the modified $2$-Laplacian defined
in Lemma \ref{lemLrad}\,.
\end{lemma}
\begin{proof}
By using the assumption $p \geq 2$
together with the comparison constellation
assumptions \eqref{eqRadConv} we obtain from
\eqref{eqLap1} and \eqref{eqLambda} that
\begin{equation*}
\begin{aligned}
&\pL^{S}(f(r(x)))  \\
&\geq  F^{p-2}(x) (p-2)\left( f''(r)\Vert
\nabla^{S}(r) \Vert^{2} +
f'(r)\Hess^{S}(r)\left(U_{r}, U_{r}
\right)\right)  \\
&\quad + F^{p-2}(x)\left(f''(r)\Vert
\nabla^{S}(r) \Vert^{2} - f'(r)\eta_{w}(r)\Vert
\nabla^{S}(r) \Vert^{2} + m\,f'(r)\left(
\eta_{w}(r) - h(r) \right)
\right) \\
& \geq  F^{p-2}(x)(p-1)\Vert
\nabla^{S}(r)\Vert^{2}
\left( f''(r) - f'(r)\eta_{w}(r)\right) \\
&\quad  + F^{p-2}(x)
f'(r)\left((p-2+m)\eta_{w}(r) - (p-2)\lambda(r)
- m\,h(r) \right) \\
& = F^{p-2}(x) \left( \left(f''(r) -
f'(r)\eta_{w}(r)\right)(p-1)\Vert
\nabla^{S}(r)\Vert^{2} +
f'(r)\mathcal{M}(r)\right).
\end{aligned}
\end{equation*}
Since $f(r)$ satisfies inequality \eqref{parenth},
we have, via $\Vert \nabla^{S}(r)\Vert \leq
g(r)$, that:
\begin{equation*}
\begin{aligned}
&\pL^{S}(f(r(x))) \\
& \geq  F^{p-2}(x)\left( \left(f''(r) -
f'(r)\eta_{w}(r)\right)(p-1)g^{2}(r) +
f'(r)\mathcal{M}(r) \right)\\
&=(p-1)F^{p-2}(x)g^{2}(r)\left( f''(r) -
f'(r)\eta_{w}(r) +
f'(r)\frac{\mathcal{M}(r)}{(p-1)g^{2}(r)}
\right) \\
&=(p-1)F^{p-2}(x)g^{2}(r)\left(f''(r) +
f'(r)\left( \frac{\mathcal{M}(r)}{(p-1)g^{2}(r)}
-
\eta_{w}(r)\right) \right) \\
&= (p-1)F^{p-2}(x)g^{2}(r) \Lmod(f(r)),
\end{aligned}
\end{equation*}
as claimed in the lemma.
\end{proof}

\section{First proof of Theorem \ref{thmMain}}\label{firstProofmain}

Next we show that (\ref{eqTransienceCond}) is
also a sufficient condition for $p$-hyperbolicity
of $S^{m}$. First we transplant the model space
solutions $\psi_{\rho,R}(r)$ of equation
(\ref{eqDir}) into the extrinsic annulus
$A_{\rho,R}=D_{R}(o)\setminus \bar D_{\rho}(o)$ in 
$S$ by defining
\[
\Psi_{\rho,R}\colon A_{\rho,R} \to \erre, \quad
\Psi_{\rho,R}(x)=\psi_{\rho,R}(r(x)).
\]
Here the extrinsic ball $D_{\rho}(o)$
is as in Definition \ref{defConstellat} and $D_{R}(o)$ is that 
component of $B_{R}(o)\cap S$ which contains $D_{\rho}(o)$.
Next we extend $\Psi_{\rho,R}$ to $S\cap \bar B_{\rho}(o)$
by setting 
$\Psi_{\rho,R}(x)=0$ for $x\in S\cap\bar B_{\rho}(o)$.

Using $w'(r) = \eta_{w}(r)w(r)$ and the balance
condition (\ref{eqBalance}) it is straightforward to check that
\begin{equation*}
\psi''_{\rho,R}(r) - \psi'_{\rho,R}(r)\eta_{w}(r) \leq 0 .
\end{equation*}

Since $\psi'_{\rho,R}(r) \geq 0$ and $\Lmod\psi_{\rho,R}=0$
in $A_{\rho, R}^{w}$, we obtain from Lemma \ref{lemPLapComp}
that
\begin{equation*}
\pLP\Psi_{\rho,R} \geq 0 
\quad\text{in}\quad D_{R}(o)\setminus\bar B_{\rho}(o).
\end{equation*}
Thus $\Psi_{\rho,R}$ is a $p$-subsolution in
$D_{R}(o)\setminus\bar B_{\rho}(o)$. In fact, $\Psi_{\rho,R}$
is a $p$-subsolution in the whole extrinsic ball $D_{R}(o)$
since $\Psi_{\rho,R}(x)=0$ for 
$x\in S\cap\bar B_{\rho}(o)$; see \cite[Theorem 7.25, Lemma 7.28]{HKM}.
Furthermore, for fixed $\rho$ and fixed $x\in S$,  
$\Psi_{\rho,R}(x)$ is defined for sufficiently large $R$
and it is decreasing as a function of $R$, see equation
(\ref{eqPsi}). Hence the limit function
\[
\Psi_{\rho}:=\lim_{R\to\infty}\Psi_{\rho,R}
\]
exists in $S$ and, moreover, it is positive in $S\setminus\bar B_{\rho}(o)$ by
\eqref{eqTransienceCond}. By \cite[Theorem 3.75]{HKM},
$\Psi_{\rho}$ is a $p$-subsolution in $S$. Hence $1-\Psi_{\rho}$
is a non-negative, non-constant $p$-supersolution in $S$, and therefore
$S$ is $p$-hyperbolic. This proves Theorem \ref{thmMain}.

\section{Proof of Corollaries}\label{secProofCor}

\begin{proof}[{\bf{Proof of Corollary \ref{corHadamard}}}] The balance
condition (\ref{eqBalance}) is clearly satisfied
by \eqref{eqCorHad}. Thus we only need to check
the $p$-hyperbolicity condition
(\ref{eqTransienceCond}). Since $g(r)\le 1$, we
have
\begin{equation*}
\frac{\mathcal{M}(r)}{(p-1)g^{2}(r)} >
(1+c)\sqrt{-b}
\end{equation*}
for some positive constant $c$ by
\eqref{eqCorHad}. Hence
\begin{equation*}
\Lambda(r)  \le
\frac{\sinh(\sqrt{-b}\,r)}{\sqrt{-b}}\exp\left(-\int_{1}^{r}(1+c)\sqrt{-b}\,dt
\right)
\end{equation*}
and therefore it is straightforward to check that
\begin{equation*}
\lim_{R\to\infty}\int_{\rho}^{R} \Lambda(t) \, dt
<  \infty ,
\end{equation*}
which concludes the proof.
\end{proof}

\begin{proof}[{\bf{Proof of Corollary \ref{corTroy}}}]
The assumptions amount to $g(r)\equiv 1$,
$h(r)\equiv 0$, and $\lambda(r)\equiv 0$ and the
only 'free' function is $w(r)$. In this intrinsic
setting we therefore have
\begin{equation*}
\mathcal{M}(r) = (m+p-2)\eta_{w}(r),
\end{equation*}
so that with $g(r)=1$ we get
\begin{equation*}
\begin{aligned}
\int_{\rho}^{r}\frac{\mathcal{M}(t)}{(p-1)g^{2}(t)}\,
dt &=
\tfrac{m+p-2}{p-1} \int_{\rho}^{r}
\frac{w'(t)\, dt}{w(t)}
\\
&=
\tfrac{m+p-2}{p-1} \log\frac{w(r)}{w(\rho)},
\end{aligned}
\end{equation*}
and hence
\begin{equation*}
\begin{aligned}
\Lambda(r)  &=  w(r) \exp\left(-
\tfrac{m+p-2}{p-1}\log\frac{w(r)}{w(\rho)}\right)\\
&= w(r)^{1- \frac{m+p-2}{p-1}}w(\rho)^{-
\frac{m+p-2}{p-1}} \\
&= w(r)^{-\frac{m-1}{p-1}}\,c(\rho),
\end{aligned}
\end{equation*}
where $c(\rho)$ is a constant depending on the
fixed inner radius of the annuli used in the proof of the
$p$-hyperbolicity. Then $\Lambda(r)$ has
bounded integral precisely if
\begin{equation*}
\lim_{R\to\infty}\int_{\rho}^{R}\frac{1}{w(r)^{\frac{m-1}{p-1}}}\,
dt  <  \infty,
\end{equation*}
as claimed.
\end{proof}

\section{$p$-capacity bounds} \label{secCapBounds}

In this section we give lower bounds on the
$p$-capacity of closed (compact) extrinsic balls
relative to $S^{m}$. Let $G\subset S^{m}$ be a
precompact open set such that $\bar
D_{\rho}(o)\subset G$. We recall from the
introduction that the $p$-capacity of $\bar
D_{\rho}(o)$ relative to $G$ is defined by
\[
\pCap (\bar D_{\rho}(o),G)=\inf_{v}\int_{G}\Vert \nabla^{S} v\Vert^{p}\,d\mu,
\]
where the infimum is taken over all real-valued functions
$v\in C^{\infty}_{0}(G)$, with $v\ge 1$ in $\bar D_{\rho}(o)$.
If $\partial G$ is regular for the Dirichlet problem for $p$-harmonic
functions, then there exists a unique function $u\in C(\bar G)$
which is $p$-harmonic in $G\setminus \bar D_{\rho}(o)$ such that
$u=0$ in $\bar D_{\rho}(o)$, $u=1$ in $\partial G$, and that
\[
\pCap (\bar D_{\rho}(o),G)=\int_{G}\Vert \nabla^{S} u\Vert^{p}\,d\mu.
\]
We refer to \cite[Chapter 6]{HKM} for the boundary regularity. For our purposes it is
enough to know that every open set can be exhausted by open sets with regular boundaries.

Since $u$ is $p$-harmonic in $G\setminus \bar D_{\rho}(o)$,
we have
\begin{equation}\label{eqCap1}
\pCap (\bar D_{\rho}(o),G)=
\int_{G}\langle\Vert \nabla^{S} u\Vert^{p-2}\nabla^{S} u,
\nabla^{S} \varphi\rangle\,d\mu
\end{equation}
for every function $\varphi\in W^{1,p}(G)$ which is continuous in
$\bar G$ with values $\varphi=0$ in $\bar D_{\rho}(o)$ and
$\varphi=1$ in $\partial G$.
In particular, \eqref{eqCap1} holds for all $0\le t<s\le 1$ with
the function
\[
\varphi(x)=\begin{cases}
0 & \text{if $u(x)\le t$}\\
\frac{u(x)-t}{s-t} & \text{if $t<u(x)<s$}\\
1 & \text{if $u(x)\ge s$}.
\end{cases}
\]
Applying the co-area formula
(\cite{MR1390760}, \cite[3.2.12, 3.2.46]{Fe},
\cite{Zi}) we obtain
\[
\pCap (\bar D_{\rho}(o),G)=
\frac 1{s-t}\int_{t}^{s}\left(\int_{u^{-1}(\tau)}\Vert\nabla^{S}u\Vert^{p-1}\,
d\mathcal{H}^{m-1}\right)\,d\tau.
\]
Letting $s\to t$ we finally get
\begin{equation}\label{eqCap2}
\pCap (\bar D_{\rho}(o),G)=
\int_{u^{-1}(t)}\Vert\nabla^{S}u\Vert^{p-1}\,
d\mathcal{H}^{m-1}
\end{equation}
for a.e. $t\in [0,1]$. We will use the equation
\eqref{eqCap2} to get lower bounds on the
$p$-capacity $\pCap (\bar D_{\rho}(o),D_{R}(o))$
in terms of the corresponding drifted
$2$-capacity in the model space.\\

Our main comparison estimate for the $p$-capacity
now reads  as follows:
\begin{theorem} \label{thmMainFinite}
Let $\{N^{n}, S^{m}, M_{w}^{m} \}$ denote a comparison constellation
on $[0, R],\ R>\rho$, in the sense of Definition \ref{defConstellat}.
Then
\begin{equation}\label{eq3.1}
\pCap\bigl(\bar D_{\rho}(o),D_{R}(o)\bigr) \geq
\left(\frac{\C_{\Lmod}(A_{\rho, R}^{w})}
{\Vol(\partial D_{\rho}^{w})}\right)^{p-1}
\int_{\partial D_{\rho}}\Vert \nabla^S r\Vert^{p-1} \, d\mathcal{H}^{m-1}.
\end{equation}
\end{theorem}

\begin{proof}
Let $G\subset D_{R}(o)$ be a precompact open set with regular boundary
such that $\bar D_{\rho}(o)\subset G$. Let $u\in C(\bar G)$ be $p$-harmonic
in $G\setminus\bar D_{\rho}(o)$ with $u=0$ in $\bar D_{\rho}(o)$
and $u=1$ in $\partial G$. Furthermore, let $\Psi_{\rho,R}$ be the $p$-subsolution
in $D_{R}(o)$ defined in Section \ref{firstProofmain}. By the comparison principle,
\[
u(x)\ge\Psi_{\rho,R}(x)
\]
for all $x\in D_{R}(o)$. Since $\nabla^{S} u$ is H{\"o}lder-continuous up
to the boundary $\partial D_{\rho}(o)$ by \cite{lieb} and
$u(x)=\Psi_{\rho,R}(x)=0$ for all
$x\in \bar D_{\rho}(o)$, we obtain
\begin{equation}\label{eqGradbound}
\Vert\nabla^{S}u(x)\Vert \ge \Vert\nabla^{S}\Psi_{\rho,R}(x)\Vert
\end{equation}
for all $x\in \partial D_{\rho}(o)$.
Combining \eqref{eqCap2} and \eqref{eqGradbound}, we arrive at
\begin{align*}
\pCap(\bar D_{\rho}(o),G)& \ge
\int_{\partial D_{\rho}}\Vert \nabla^S \Psi_{\rho,R}\Vert^{p-1}\, d\mathcal{H}^{n-1}\\
& =\left(\psi'_{\rho,R}(\rho)\right)^{p-1}
\int_{\partial D_{\rho}}\Vert \nabla^S r\Vert^{p-1}\, d\mathcal{H}^{m-1} \\
& = \left(\frac{\C_{\Lmod}(A_{\rho, R}^{w})}
{\Vol(\partial D_{\rho}^{w})}\right)^{p-1}
\int_{\partial D_{\rho}}\Vert \nabla^S r\Vert^{p-1}\, d\mathcal{H}^{m-1}.
\end{align*}
The desired estimate \eqref{eq3.1} now follows since
\[
\pCap\bigl(\bar D_{\rho}(o),D_{R}(o)\bigr)
=\inf_{G}\pCap(\bar D_{\rho}(o),G),
\]
where $G\subset D_{R}(o)$ is a precompact open set with regular boundary.
\end{proof}

\subsection{Second Proof of Theorem \ref{thmMain}}
\label{subsecSecondProof} Using the explicit
capacity comparison obtained in Theorem \ref{thmMainFinite}
we finally observe the following direct
proof of the main theorem.

Let $\{N^{n}, S^{m}, M_{w}^{m} \}$ denote a
comparison constellation on $[0, \infty]$ in the
sense of Definition \ref{defConstellat}. By
assumption $D_{\rho}(o)$ is precompact with a
smooth boundary and thence, in equation
(\ref{eq3.1}) we have
\[
\int_{\partial D_{\rho}}\Vert \nabla^S
r\Vert^{p-1} \, d\mathcal{H}^{m-1} \, > 0 .
\]

From  (\ref{eqModelCap}) and the assumption
(\ref{eqTransienceCond}) we also have
\[
\lim_{R \to \infty}\C_{\Lmod}(A_{\rho, R}^{w}) \,
> \, 0 ,
\]
so that Theorem \ref{thmMainFinite} implies:
\[
\pCap\bigl(\bar D_{\rho}(o), S^{m}\bigr)\, = \,
\lim_{R\to \infty} \pCap\bigl(\bar
D_{\rho}(o),D_{R}(o)\bigr)\, > \, 0 .
\]

Thus $\bar D_{\rho}(o)$ is a compact subset with
positive $p$-capacity in $S^{m}$, and
$p$-hyperbolicity of that submanifold follows
again.

\enddocument